\documentstyle[12pt]{article}
\input tcilatex
\QQQ{Language}{
American English
}

\begin{document}

\title{Comments on Lagrange Partial Differential Equation}
\author{C. Viazminsky \\
%EndAName
IITAP, Iowa State University, Ames, IA 50011\\
and Dept. of Physics, University of Aleppo,\\
Syria}
\maketitle

\begin{abstract}
The relations between solutions of the three types of totally linear partial
differential equations of first order are presented. The approach is based
on factorization of a non-homogeneous first order differential operator to
products consisting of a scalar function, a homogeneous first order
differential operator and the reciprocal of the scalar function. The
factorization procedure is utilized to show that all totally linear
differential equations of first order can be transformed to each other, and
in particular to a homogeneous one.
\end{abstract}

\section{Introduction}

The method for solution of Lagrange partial differential equation is well
known, and is found almost in every text book on partial differential
equations\cite{Piaggio,Ayres,Carrier}. Our goal in this article is to show
how the factorization of a non-homogeneous first order differential operator
leads quite naturally to simple relations between the solutions of three
related types of Lagrange equations.

Let $E$ be an open subset of $\Re ^n,$ $L$ be a continuous vector field in $E
$, and denote by $C^k(E)$ the set of real valued functions which are
continuously differentiable of order $k$ on $E.$ The continuous vector field 
$L$ may be viewed as a differential operator \cite{Abraham} from $C^1(E)$ to 
$C^0(E)$. For each real valued continuous function $q:E\rightarrow R$ there
corresponds an operator $L+q:C^1(E)\rightarrow C^0(E)$ defined by $(L+q)\psi
=L\psi +q\psi ,\forall \psi \in C^1(E)$, where $L\psi $ is the Lie
derivative of the function $\psi $ with respect to the field $L,$ and $q\psi 
$ is the usual product of two functions. Our goal in this work is to study
the relations between the solutions of the partial differential equations

$\;\;(i)\;L\phi =0$ \ \ \ \ \ $(ii)\ (L+q)\psi =0,\;\;\;\;(iii)\;(L+q)\chi
=b,$ \ \ \ \ \ \ \ \ \ \ \ \ \ \ \ \ \ \ \ \ \ \ \ \ \ \ \ \ \ \ \ \ \ \ \ \
\ \ \ \ \ \ \ \ \ \ \ \ \ \ \ \ \ \ \ \ \ \ \ \ \ \ \ \ \ \ \ \ \ \ \ \ \ \
\ \ \ 

where $b$ is a continuous function on $E$. The approach followed here hinges
on factorization of the first order non-homogeneous differential operator $%
(L+q)$ to a product of a scalar function, a homogeneous differential
operator, and the reciprocal scalar function.

\section{Factorization of a First Order \protect\\Non-Homogeneous Operator}

Let $\eta \in C^1(E)$ be a non-zero solution of the differential equation 
\begin{equation}
(L+q)\eta =0.  \label{e1}
\end{equation}
Equivalently,$\,\eta $ is any element in the kernel of the linear operator $%
L+q$ that is different from zero. As a first step we assume that $\eta \,$%
has no zeros in $E$, and hence $\eta ^{-1}\,$exists and of class $C^1(E)$.\
The general case in which $\eta $ vanishes on a subset $\delta \subset E$
will be considered in section 4. We start by proving a useful operator
equality on which hinges the method of reducing one type of Lagrange
equations to another.

\begin{theorem}
In $C^1(E)$ the following operator equality holds 
\begin{equation}
\eta L\eta ^{-1}=L+q.  \label{e2}
\end{equation}
Proof: for every $\psi \in C^1(E).$
\end{theorem}

\[
(\eta L\eta ^{-1})\psi =\eta (\eta ^{-1}L-\eta ^{-2}(L\eta ))\psi =(L-\eta
^{-1}(L\eta ))\psi =(L+q)\psi . 
\]
We have used equation (\ref{e1}) to make the last step.%
\setcounter{theorem}{0}

\begin{corollary}
Equality (\ref{e3}) is equivalent to 
\begin{equation}
L=\eta ^{-1}(L+q)\eta   \label{e3}
\end{equation}
which shows that all operators of the form $(L+q)\;$which are based on the
same field $L$ may be transformed to $L,$ and accordingly to each other.\ 
\end{corollary}

\begin{corollary}
The equality (\ref{e2}) shows that the left hand-side must not be dependent
on the particular solution $\eta $ of equation (\ref{e1})$,$ since its right
hand-side is not. Therefore if $\xi $ is another solution of (\ref{e1}),
then by equation (\ref{e2}) and a similar equation written for the solution $%
\xi ,$ we have$\;\xi ^{-1}\eta L\eta ^{-1}\xi =L.$ This yields $L(\xi /\eta
)=0.$ \ \ \ \ \ \ \ \ \ \ \ \ \ \ \ \ \ \ \ \ \ \ \ \ \ \ \ \ \ \ \ \ \ \ \
\ \ \ \ \ \ \ \ \ \ \ \ \ \ \ \ \ \ \ \ \ \ \ 
\end{corollary}

\begin{corollary}
From (\ref{e3}) we deduce that 
\begin{equation}
L^k=\eta ^{-1}(L+q)^k\eta \,\,\,\,\,and\;(L+q)^k=\eta L^k\eta ^{-1}
\label{e4}
\end{equation}
$\;$where $k$ is a non-negative integer. If $L\;$is invertible then the
latter relation holds for all integers.\thinspace \thinspace We assume in
relation (\ref{e4})\thinspace \thinspace \thinspace that $L$ is a $C^{k-1}\;$%
vector\ field,\thinspace $q$ is a $C^{k-1}$ function and $\eta \;$is a $C^k$
function. \thinspace \thinspace \ \ \ \ \ \ \ \ \ \ \ \ \ \ \ \ \ \ \ \ \ \
\ \ \ \ \ \ \ \ \ \ \ \ \ \ \ \ \ \ \ \ \ \ \ \ \ \ \ \ \ \ \ \ \ \ \ \ \ \
\ \ \ 
\end{corollary}

\begin{corollary}
If (\ref{e1}) holds then it is easily checked that $(L+kq)\eta ^k=0$, and
hence 
\begin{equation}
L=\eta ^{-k}(L+kq)\eta ^k  \label{e5}
\end{equation}
$\;$In general, and for any real number $\alpha ,$ we have 
\begin{equation}
L=\mid \eta \mid ^{-\alpha }(L+\alpha q)\mid \eta \mid ^\alpha   \label{e6}
\end{equation}
\end{corollary}

\begin{corollary}
If $Q$ is a real valued continuous function on $E$ then 
\begin{equation}
\eta ^{-1}(L+q+Q)\eta =L+Q.  \label{e7}
\end{equation}
\end{corollary}

\begin{corollary}
Take $Q=-\lambda \,(\lambda \in R)$ in corollary (\ref{e4}) to obtain 
\begin{equation}
(L-\lambda )\psi _\lambda =0\Leftrightarrow (L+q-\lambda )(\eta \psi
_\lambda )=0  \label{e8}
\end{equation}
$\;$The last relation states that: if $\psi _\lambda $ is an eigenfunction
of the operator $L\;$belonging to the eigenvalue $\lambda $ then$\;\eta \psi
_\lambda $ is an eigenfunction of the operator $L+q$ belonging to the same
eigenvalue $\lambda .\;$
\end{corollary}

\begin{example}
Take $L+q=\frac d{dx}+2x:C^1(R)\rightarrow C^0(R).\;$Since $\eta =e^{-x^2}$
is a solution of (\ref{e1}), we have 
\[
\frac d{dx}+2x=e^{-x^2}\frac d{dx}e^{x^2}
\]
It is obvious that every complex number $\lambda $ is an eigenvalue of the
operator $\frac d{dx}$ to which an eigenfunction $\psi _\lambda =e^{\lambda
x}$ belongs. In accordance with the last corollary, it is easily checked
that $\lambda $ is an eigenvalue of the operator $\frac d{dx}+2x$ to which
the function $\eta \psi _\lambda =e^{-x^2+\lambda x}$ belongs.
\end{example}

\section{On Totally Linear Partial Differential Equations}

Let $(x_1,...,x_n)$ be a global system of coordinates on the region $E$, in
which $L$ is expressed as 
\begin{equation}
L=\stackrel{n}{\stackunder{k=1}{\sum }}a_k(x_1,...,x_n)\partial /\partial
x_k,  \label{e9}
\end{equation}
where the components $a_k$ are of class $C^0$ on $E$. We shall describe the
partial differential equation 
\begin{equation}
(L+q)\chi =b  \label{e10}
\end{equation}
where $q(x)$ and $b(x)$ are continuous functions on $E$, as totally linear.
The totally linear equation 
\begin{equation}
(L+q)\psi =0,  \label{e11}
\end{equation}
will be referred to as the non-homogeneous reduced equation corresponding to
(\ref{e10}), or simply, as the non-homogeneous equation. Equation (\ref{e10}%
) is a special type of Lagrange equation. The method of solution of Lagrange
equation (\ref{e10}), and consequently equations (\ref{e10}) and (\ref{e11})
is well known \cite{Piaggio}. However we aim here to utilize equality (\ref
{e2}) to reduce the non-homogeneous equation (\ref{e11}) to the homogeneous
equation 
\begin{equation}
L\phi =0,  \label{e12}
\end{equation}
and express its general solution in terms of a particular solution and the
general solution of (\ref{e12}). Alternatively, to express the general
solution of (\ref{e11}) in terms of $n$ particular solutions. The results we
have just pointed to are expressed in the following facts in which we assume
that $\eta $ is a solution of (\ref{e11}) on $E$ and that it has no zeros on 
$E$.

F1. A function $\phi _0$ is a solution on $E$ of $L\phi =0$ if and only if $%
\psi _0=\eta \phi _0$ is a solution of $(L+q)\psi =0$ on $E$.

The proof is a direct consequence of corollary 1 in the previous section.

F2. Let $Q:E\rightarrow \Re $ be continuous. By corollary 5 in the previous
section, a function $\psi _0$ is a solution of (\ref{e11}) on $E$ iff $\Psi
_0=\eta \psi _0$ is a solution of $(L+q+Q)\Psi =0$ on $E.$ In a more
familiar language to the subject of differential equation, the
transformation $\Psi =\eta \psi $ reduces the last equation to (\ref{e11}).

F3. If $\eta $ and $\xi $ are solutions of (\ref{e11}) then $\xi /\eta $ is
a solution of (\ref{e12}). Indeed, from equation (\ref{e11}) which is
satisfied by $\eta $ and $\xi $ we get $\eta L\xi =\xi L\eta $, and hence 
\[
L(\xi /\eta )=\eta ^{-2}(\eta L\xi -\xi L\eta )=0.
\]

F4. The general solution of the reduced non-homogenous equation (\ref{e11})
is given by 
\begin{equation}
\psi =\eta \;f(\phi _1,....,\phi _{n-1}),  \label{e13}
\end{equation}
$\;$where 
\begin{equation}
\phi _i=\phi _i(x_1,...,x_n)\;\;\;\;\;\;\;(i=1,....,n-1)  \label{e14}
\end{equation}
are $(n-1)$ functionally independent solutions of (\ref{e12}).

Proof. According to the standard method in solving Lagrange equation \cite
{Piaggio}, the general integral of the homogeneous equation (\ref{e12}) is
given by 
\begin{equation}
\phi =f(\phi _1,....,\phi _{n-1}),  \label{e15}
\end{equation}
where $f$ is an arbitrary $C^1$function in its arguments. Now if $\psi $ is
a solution of (\ref{e11}) then by F3 $\psi /\eta $ is a solution of (\ref
{e12}), and hence it must be of the form (\ref{e15}). It follows that the
general solution of (\ref{e11}) is given by (\ref{e13}).

F5. If 
\begin{equation}
\eta (x_1,....,x_n),\;\eta _j(x_1,....,x_n)\;\;\;\;\;\;(j=1,....,n-1)
\label{e16}
\end{equation}
are functionally independent solutions of (\ref{e11}) then the ratios 
\begin{equation}
\phi _j=\eta _j/\eta \;\;\;\;(j=1,....,n-1)  \label{e17}
\end{equation}
are solutions of (\ref{e12}). It is easily seen that these ratios are
functionally independent. Hence 
\begin{equation}
\phi =f(\eta _1/\eta ,....,\eta _{n-1}/\eta )  \label{e18}
\end{equation}
is the general solution of the homogeneous equation (\ref{e12}), and 
\begin{equation}
\psi =\eta \;f(\eta _1/\eta ,....,\eta _{n-1}/\eta )  \label{e19}
\end{equation}
is the general solution of the non-homogeneous equation (\ref{e11}). If $%
\eta _n$ is a further solution of (\ref{e11}) then by (\ref{e19}) 

\begin{equation}
\eta _n/\eta =f_0(\eta _1/\eta ,....,\eta _{n-1}/\eta ),  \label{e20}
\end{equation}
where $f_0$ is some specified $C^1$ function.

F6. Let $\chi _0$ be a solution of the totally linear equation (\ref{e10}).
If $\chi $ is any other solution of (\ref{e10}) then $(L+q)(\chi -\chi _0)=0.
$ Hence $\psi =\chi -\chi _0$ is a solution of (\ref{e11}), and consequently
must be of the form (\ref{e13}). It follows that the general solution of (%
\ref{e10}) is given by 
\begin{equation}
\chi =\chi _0+\eta \;f(\phi _1,....,\phi _{n-1}).  \label{e21}
\end{equation}
If a second particular solution $\chi _1$ of (\ref{e10}) is given, then $%
\eta =\chi _1-\chi _0$ is a solution of (\ref{e11}), and therefore the
general solution of (\ref{e10}) can be written as follows 
\begin{equation}
(\chi -\chi _0)/(\chi _1-\chi _0)=f(\phi _1,....,\phi _{n-1}).  \label{e22}
\end{equation}
If $\eta _2$ is a third solution of (\ref{e10}), then by (\ref{e22}), ($\chi
_2-\chi _1)/(\chi _1-\chi _0)$ is a solution of (\ref{e12}) and has
accordingly the form (\ref{e15}). If it is known $(n+1)$ solutions $\chi
_0,\chi _1,.....,\chi _n$ of (\ref{e10}), then 
\[
\eta =\chi _{1-}\chi _0,\;\eta _1=\chi _2-\chi _0,....,\;\eta _{n-1}=\chi
_n-\chi _0
\]
are solutions of (\ref{e11}), and hence 
\begin{equation}
\phi _1=\frac{\eta _1}\eta =\frac{\chi _2-\chi _0}{\chi _1-\chi _0}%
\;,......,\phi _{n-1}=\frac{\eta _{n-1}}\eta =\frac{\chi _n-\chi _0}{\chi
_1-\chi _0}  \label{e23}
\end{equation}
$\;$are solutions of (\ref{e12}). If these ratios are functionally
independent then 
\begin{equation}
\frac{\chi -\chi _0}{\chi _1-\chi _0}=f(\frac{\chi _2-\chi _0}{\chi _1-\chi
_0}\;,.....,\;\frac{\chi _n-\chi _0}{\chi _1-\chi _0})  \label{e24}
\end{equation}
$\;$gives the general solution $\chi $ of (\ref{e10}) in terms of $(n+1)$
particular solutions. It is clear of course that these $(n+1)$ solutions of (%
\ref{e10}) are certainly functionally dependent. If $\chi _{n+1}$ is a
further solution of (\ref{e10}) then by (\ref{e24}) $(\chi _{n+1}-\chi
_0)/(\chi _1-\chi _0)$ is some function of the solutions (\ref{e23}) of
equation (\ref{e12}). 

\section{The General Case\ }

We proceed here to study the general case in which the solution $\eta $ of
equation (\ref{e1}) is defined only on an open subset $E_\eta \subset E,$
and vanishes on a closed subset $\delta \subset E$. This latter assumption
embodies most interesting cases which are encountered in practical examples,
such as $\eta $ vanishing at a finite or a countable set, or outside an open
subset of its domain of definition. In the present case the operator $\eta
L\eta ^{-1}$ exists only on $E_\eta -\delta .$ Hence the equality (\ref{e2})
must be replaced by the inclusion relation $(\eta L\eta ^{-1}\subset L+q);$
it is an equality only on $E_\eta -\delta .$\ 

Some minor modification has to be made so that the corollaries of theorem1
in section 2 and the facts in section 3 conform to the new assumptions. As
an example the fact F1 in section 3 must be modified as follows:

F1. Let $\eta $ be a solution of (\ref{e1}) on $E_\eta $ which vanishes on a
closed subset $\delta \subset E_\eta $.

(i) $\phi $ is a solution of (\ref{e12}) on $E_\phi $ $\Rightarrow \eta \phi 
$ is a solution of (\ref{e11}) on $E_\phi \cap E_\eta .$

(ii) $\psi $ is a solution of (\ref{e11}) on $E_\psi $ $\Rightarrow \psi
/\eta $ is a solution of (\ref{e12}) on $E_\psi \cap E_\eta -\delta .$

Proof: (i) on $E_\phi \cap E_\eta $ $,$ where $\eta $ is defined, we have 
\[
(L+q)(\eta \phi )=\phi (L+q)\eta +\eta L\phi .
\]

If $L\phi =0$ on $E_\phi $ then $(L+q)(\eta \phi )=0$ on $E_\eta \cap E_\phi
.$

(ii) On $E_\phi \cap E_\eta -\delta ,$ where $\psi /\eta $ is defined and is
of class $C^1,$ we have $L(\psi /\eta )=0$, as it was shown in F3.

It must be noted that the last fact determines the smallest domains on which 
$\phi $ and $\psi $ are defined. These domains may be extended by continuity
to larger domains.

The remaining facts and corollaries can be modified in a similar way.

Example. Consider the operator $L=x\partial /\partial x+y\partial /\partial
y+z\partial /\partial z$ which is continuous on $\Re ^3$. The general
solution of the equation $L\phi =0$ is $\phi =f(y/x,z/x),$ where $f$ is an
arbitrary $C^1$ function. The equation $(L+3/2)\psi =0$ admits the
particular solution $\eta =x^{-3/2}$ on $E_\eta =\{(x,y,z)\in \Re ^3:x>0\}.$
The function $\phi =x/y$ is a solution of $L\phi =0$ on $E_\phi
=\{(x,y,z)\in \Re ^3:y\neq 0\},$ and the function $\psi =\eta \phi
=x^{-1/2}y^{-1}$ is a solution of $(L+3/2)\psi =0$ on $\{(x,y,z)\in \Re
^3:y\neq 0,x>0\}=E_\eta \cap E_\phi .$ It is clear that we could have
adopted $\eta ^{\prime }=\left| x\right| ^{-3/2},$ and $\phi ^{\prime
}=\left| x\right| /y$ (since $L$ is linear), to obtain a solution $\psi
^{\prime }=\left| x\right| ^{-1/2}y^{-1}$  defined on $\Re
^3-\{x=0,y=0\}=E_{\eta ^{\prime }}\cap E_{\phi ^{^{\prime }}}.$ If we
replace in the given equation $3/2$ by $1$ , then $\eta =1/x,$ and hence $%
\eta \phi =1/y$ is a solution of $(L+1)\psi =0$ on $E_\phi .$ 

Consider now the totally linear equation $(L+1)\chi =3x^2.$ Four particular
solutions of this equation are $\chi _0=x^2$ , $\chi _1=x^2+x^{-1},\;\chi
_2=x^2+y^{-1},\;\chi _3=x^2+z^{-1}.$ It can be checked that the formula (\ref
{e24}) yields, in accordance with (\ref{e21}), the general solution $\chi
=x^2+x^{-1}f(y/x,z/x).$

\section{An Algebraic Comment}

Let $A$ denote the commutative real algebra formed by the set of all
functions defined in $E$. The set of solutions of the homogeneous equation $%
L\phi =0,$ namely $KerL,$ is clearly a subalgebra of $A$. The set of
solutions of non-homogeneous equation (\ref{e11}) is a coset of this
subalgebra determined by a particular solution of (\ref{e11}). If $\eta $ is
a particular solution of (\ref{e11}) defined every where in $E$, then $%
Ker(L+q)=\eta \;KerL.$ It is clear that $Ker(L+q)$ is a sub-vector space of $%
A$. The set of solutions of equation (\ref{e10}) is a coset of the
sub-vector space $Ker(L+q)$ determined by any particular solution of (\ref
{e10}).


\begin{thebibliography}{9}
\bibitem{Abraham}  Abraham R and Marsden J, Reading, Mass.,
Benjamin/Cummings Pub. Co. (1978).

\bibitem{Piaggio}  Piaggio H, An Elementary Treatise on Differential
Equations and Their Applications. G. Bell and Sons LTD, London (1969).

\bibitem{Ayres}  Ayres F. Differential Equations. Schaum's Outline Series.
McGraw Hill Book Company, New York (1952).

\bibitem{Carrier}  Carrier G F, Partial Differential Equations. Academic
Press, New York (1979).
\end{thebibliography}
\end{document}